%% file: cramerwold-arxiv.tex
\def\versiondate{22 April 2017}
\input math10.macros
\magnification=1140

\let\nobibtex = t
\let\noarrow = t
\input eplain
\beginpackages
 \usepackage{url}
 \usepackage{color}  
 \usepackage{graphicx}
\endpackages

\input Ref.macros

\checkdefinedreferencetrue
\continuousfigurenumberingtrue
\theoremcountingtrue
\sectionnumberstrue
\forwardreferencetrue
\citationgenerationtrue
\nobracketcittrue
\hyperstrue
\initialeqmacro


\input cramerwold.key
\bibsty{../../texstuff/myapalike}

\def\bfz{{\bf 0}}

\def\ip#1,#2{\langle #1, #2 \rangle}
\def\Bigip#1,#2{\Big\langle #1, #2 \Big\rangle}
\def\bigip#1,#2{\big\langle #1, #2 \big\rangle}

\def\sphere{\Bbb S}
\def\sfc{\Omega} 
\def\hfsp{{\scr S}}  
\def\csp{{\rm c}}  
\def\n{{\bf n}}
\def\bfz{{\bf 0}}

\ifproofmode \relax \else\head{}
{Version of \versiondate}\fi 
\vglue20pt

\title{A Calculus Proof of the Cram\'er--Wold Theorem}

\author{Russell Lyons and Kevin Zumbrun}

\abstract{We present a short, elementary proof not involving Fourier transforms
of the theorem of Cram\'er and Wold that a Borel probability measure is
determined by its values on half-spaces.
}

\bottomII{Primary 
60E10. 
Secondary
44A12, 
53C65. 
}
{Distributions of probability measures, Radon transform, integral geometry.}
{Research partially supported by NSF grants DMS-1007244 and DMS-1400555.}

\bsection{Introduction}{s.intro}

In this note, we give a brief and elementary proof, not involving Fourier
transforms, of a theorem of Cram\'er and Wold.

The fundamental theorem of \ref b.CramerWold/ states that a Borel probability
measure on Euclidean space is determined by the values it
assigns to all half-spaces (equivalently, by its projections to lines through
the origin).
This theorem is proved easily with the aid of Fourier analysis. However,
generations of probabilists have learned from the editions of the
textbook of \ref b.Billingsley:PM/ that despite the elementary statement of
the theorem,
no proof was known that did not use Fourier transforms.
That changed with the publication of
\ref b.walther:withAdd/, who used
Gaussians.
Walther's proof depends on
a nice idea, but its implementation uses $1 {1 \over 2}$ pages of
calculations. See Section 8.7 of \ref b.Pollard/ for another presentation
of Walther's proof.
By contrast, our proof uses only natural constructions and avoids
calculations.

A brief and somewhat inaccurate outline of our proof is the following.
Using Crofton's measure on half-spaces, we show that knowledge of $\mu(S)$
for all half-spaces $S$ determines the $\mu$-average distance $f_\mu(x)$
to every point, $x$.
We then show that a suitable power of the Laplacian applied to $f_\mu$
yields a constant times $\mu$.
Thus, integral geometry combined with differentiation recovers $\mu$.

\bsection{Proof}{s.proof}

Let $\hfsp$ be the set of closed half-spaces $S \subset\R^n$.

\proclaim The Cram\'er--Wold Theorem. Let $\mu$ and $\nu$ be Borel
probability measures on $\R^n$ such that $\mu(S) = \nu(S)$ for all $S \in
\hfsp$. Then $\mu = \nu$.

\proof
Let
$\sigma$ be the (infinite) Borel measure on $\hfsp$ 
that is invariant under isometries, normalized so that
$$
\sigma\big(\{0 \in S, x \notin S\}\big) = \|x\|/2
\label e.distance
$$
for $\|x\|= 1$.
The measure $\sigma$ goes back to \ref b.Crofton/ (in two dimensions); it
can be constructed as follows. (See Theorem 5.1.1 of \ref b.SW:book/ for
a generalization.)
Let $\sfc_{n-1}$ denote hypersurface area measure on the unit sphere
$\sphere^{n-1} \subset \R^n$, and let
$\lambda$ denote Lebesgue measure on $\R$.
Write $\varphi \colon \sphere^{n-1} \times \R \to \hfsp$ for the map
$$
\varphi(\omega, p)
:=
\bigl\{ x \in \R^n \st \ip \omega, x \ge p \bigr\}
\,.
$$
Then $\sigma := \alpha_n \cdot \varphi_*(\sfc_{n-1} \times \lambda)$ for some
constant $\alpha_n$ whose value does not concern us.
It is clear that $\sigma$ is invariant under rotations about the origin
and under reflections in hyperplanes that pass through the origin.
Translation invariance
amounts to the property that for $y \in \R^n$, the pushforward by
$\varphi_y(\omega, p) := \varphi(\omega, p) - y$ is the same measure. But
since 
$$\eqaln{
\varphi(\omega, p) - y 
&= 
\bigl\{ x - y \in \R^n \st \ip \omega, x \ge p \bigr\}
=
\bigl\{ x \in \R^n \st \ip \omega, {x + y} \ge p \bigr\}
\cr&=
\bigl\{ x \in \R^n \st \ip \omega, x \ge p - \ip \omega, y \bigr\}
=
\varphi\bigl(\omega, p - \ip \omega, y\bigr)
\,,
}$$
isometry invariance of $\lambda$ gives this property.
The isometry invariance of $\sigma$ implies that $\sigma\big(\{0 \in S,\, x
\notin S\}\big)$ is a function of $\|x\|$ alone; additivity for collinear
segments shows that it is a linear function. Thus, we may choose $\alpha_n$
so that \ref e.distance/ holds.

From \ref e.distance/ and isometry invariance, we have
$$
\|x\| =
\int_\hfsp |1_S(0) - 1_S(x)|^2 \,d\sigma(S)
\,.
$$
Integrating with respect to a
signed measure $\mu$ on $\R^n$ with compact support, we obtain
$$
\int_{\R^n} \|x\| \,d\mu(x)
=
\int_\hfsp \int_{\R^n} |1_S(0) - 1_S(x)|^2 \,d\mu(x) \,d\sigma(S)
=
\int_\hfsp \Bigl[1_S(0)\big(1 - 2\mu(S)\big) + \mu(S)\Bigr] \,d\sigma(S)
\,.
$$
The choice of 0 was arbitrary, so making another
choice and subtracting, we get
$$
\int_{\R^n} \big(\|y - x\| - \|x\|\big) \,d\mu(x)
=
\int_\hfsp \Bigl[\bigl(1_S(y) - 1_S(0)\bigr)\bigl(1 - 2\mu(S)\bigr)\Bigr]
\,d\sigma(S)
\,.
$$
By taking a limit, we see that this equation holds for every finite signed
measure, $\mu$.

Define
$$
f_\mu(y)
:=
\int_{\R^n} \big(\|y - x\| - \|x\|\big) \,d\mu(x)
\,.
$$
We have shown that the function $S \mapsto \mu(S)$ determines $f_\mu$. It
remains to show that $f_\mu$ determines $\mu$.

The idea is that if $n = 2m-1$ is odd, then $\Delta^m f_\mu = c_m
\mu$ for some constant $c_m$, using the fundamental solution of the
Laplacian, $\Delta$. This then establishes the Cram\'er--Wold theorem
in odd dimensions.
But since an even dimension embeds in the next higher dimension, 
the Cram\'er--Wold theorem follows in even dimensions as well.
That is, we may identify a measure $\mu$ on $\R^{2m}$ with a measure $\mu'$
on $\R^{2m} \times \{0\} \subset \R^{2m+1}$. 
The function $S \mapsto
\mu(S)$ on half-spaces $S \subset \R^{2m}$ determines the values
$\mu'(S')$ for half-spaces $S' \subset \R^{2m+1}$. Since this determines
$\mu'$, the theorem follows for $\mu$.

We now show that
$\Delta^m f_\mu = c_m \mu$ in an appropriate sense for $\mu$ on $\R^{2m-1}$.
Recall Green's second identity, which says that for a
bounded domain $D \subset \R^n$ with $C^1$ boundary $\partial D$ having
outward unit normal $\n$ and two functions $\phi, \psi \in C^2(\overlinesl
D\,)$, we
have
$$
\int_D (\phi \Delta \psi - \psi \Delta \phi) 
=
\int_{\partial D} (\phi \nabla_\n \psi - \psi \nabla_\n \phi)
\,.
$$
Recall also that if $F \colon \R^n \to \R$ is such that $F(x) =
G\bigl(\|x\|\bigr)$ depends only on $r := \|x\|$, then 
$$
(\Delta F)(x)
=
G''(r) + {n-1 \over r} G'(r)
\,.
$$
In particular, $\Delta r^k = k (k + n - 2) r^{k-2}$.
If the support of $\psi$ lies in the interior of a ball $B(\bfz, R)$ 
and $\phi(x) = r^k$ with $k > -n + 2$, then letting $D$ be $B(\bfz,
R) \setminus B(\bfz, \epsilon)$ with $\epsilon \to 0$ shows that
$\int_{\R^n} \phi \Delta \psi = \int_{\R^n} \psi \Delta \phi$. Similarly,
if $k = - n +2$, then 
$\int_{\R^n} \phi \Delta \psi = \beta_{n-1} \psi(\bfz)$, where
$\beta_{n-1}$ is the surface area of $\sphere^{n-1}$.

To show that $f_\mu$
determines $\int g \,d\mu$ for all $g \in C^\infty_{\csp}(\R^{2m-1})$, we
now prove that with $c_m := 2(-2\pi)^{m-1} (2m-2)!!$, where $!!$ denotes
the double factorial, we have
$$
\int g \,d\mu = c_m^{-1}
\int_{\R^{2m-1}} f_\mu(y) (\Delta^m g)(y) \,d\lambda(y)
\,,
$$
where now
$\lambda$ denotes Lebesgue measure on $\R^{2m-1}$. Fubini's theorem
yields
$$
\int_{\R^{2m-1}} f_\mu(y) (\Delta^m g)(y) \,d\lambda(y) 
=
\int_{\R^{2m-1}} \int_{\R^{2m-1}} \big(\|y - x\| - \|x\|\big) (\Delta^m g)(y)
\, d\lambda(y) \,d\mu(x)
\,.
$$
Applying the preceding Green formulas (translated to $x$)
repeatedly to the inner integral, we obtain
$$\eqaln{
\int_{\R^{2m-1}} \big(\|y - x\| - \|x\|\big) (\Delta^m g)(y)
\, d\lambda(y) 
&=
\int_{\R^{2m-1}} \Delta_y^{m-1}
\big(\|y - x\| - \|x\|\big) \Delta g(y) \,
d\lambda(y) 
\cr&=
c_m g(x) 
\,,
}$$
as desired.
\Qed

\eject

Our inversion formula $\mu = c_m^{-1} \Delta^m f_\mu$ in $\R^{2m-1}$
is similar to a well-known inversion formula for the Radon
transform due to \ref b.Radon/ and \ref b.John:book/, p.~13:
If $f \in C^1_{\csp}(\R^n)$, then writing $J(\omega, p) := \int_{\langle
\omega, x \rangle = p} f(x) \,dx$ for the integral of $f$ on a hyperplane
perpendicular to $\omega\in \sphere^{n-1}$, we have
$$
f(x) = \cases{\vadjust{\kern1pt}%
{1 \over 2} (2\pi)^{1-n} (-\Delta_x)^{(n-1)/2}
\int_{\sphere^{n-1}} J(\omega,\langle\omega, x\rangle)\,
d\Omega_{n-1}(\omega) &if $n$ is
odd\cr
- (2\pi)^{-n} (-\Delta_x)^{(n-2)/2}
\int_{\sphere^{n-1}} \int_\R {dJ(\omega, p) \over p - \langle\omega,
x\rangle}\, d\Omega_{n-1}(\omega) &if $n$ is
even.\cr}
$$
Apparently it was not realized until pointed out by \ref b.Renyi/ that the
theorem of \ref b.CramerWold/ generalized the injectivity
results of Radon, John, and others.

The injectivity of the map $\mu \mapsto \int_{\R^n} \|x\| \,d\mu(x)$ for
probability measures $\mu$ with finite first moment holds in other spaces
as well. On metric spaces of negative type, it is equivalent to strong
negative type. See Remark 3.4 of \ref b.L:dcov/ for details and references.
That paper also shows its relevance to statistics.

\def\noop#1{\relax}
\input cramerwold.bbl

\filbreak
\begingroup
\eightpoint\sc
\parindent=0pt\baselineskip=10pt

Department of Mathematics,
831 E. 3rd St.,
Indiana University,
Bloomington, IN 47405-7106
\emailwww{rdlyons@indiana.edu}
{http://pages.iu.edu/\string~rdlyons/}
\emailwww{kzumbrun@indiana.edu}
{http://pages.iu.edu/\string~kzumbrun/}

\endgroup

\bye

%% file: cramerwold.bbl
\def\polhk#1{\setbox0=\hbox{#1}{\ooalign{\hidewidth
  \lower1.5ex\hbox{`}\hidewidth\crcr\unhbox0}}}
  \def\soft#1{\leavevmode\setbox0=\hbox{h}\dimen7=\ht0\advance \dimen7
  by-1ex\relax\if t#1\relax\rlap{\raise.6\dimen7
  \hbox{\kern.3ex\char'47}}#1\relax\else\if T#1\relax
  \rlap{\raise.5\dimen7\hbox{\kern1.3ex\char'47}}#1\relax \else\if
  d#1\relax\rlap{\raise.5\dimen7\hbox{\kern.9ex \char'47}}#1\relax\else\if
  D#1\relax\rlap{\raise.5\dimen7 \hbox{\kern1.4ex\char'47}}#1\relax\else\if
  l#1\relax \rlap{\raise.5\dimen7\hbox{\kern.4ex\char'47}}#1\relax \else\if
  L#1\relax\rlap{\raise.5\dimen7\hbox{\kern.7ex
  \char'47}}#1\relax\else\message{accent \string\soft \space #1 not
  defined!}#1\relax\fi\fi\fi\fi\fi\fi} 
  \def\lfhook#1{\setbox0=\hbox{#1}{\ooalign{\hidewidth
  \lower1.5ex\hbox{'}\hidewidth\crcr\unhbox0}}} 
   
\def\temp{\let\linkit=\linkyear \apaliketrue}
\temp
\ifcitationgeneration\immediate\write\labelfile{\sanitize\temp}\fi
\def\startreferences{
 \vskip0pt plus.3\vsize \penalty -150 \vskip0pt
 plus-.3\vsize \bigskip\bigskip \vskip \parskip
 \begingroup\baselineskip=12pt\frenchspacing
 \bibliographytitle
 \vskip12pt\parindent=0pt
 \def\and{{\rm and}}
 \def\em{\it}
 \def\newblock{\hskip .11em plus.33em minus.07em}
 \def\bibauthor##1{{\sc ##1}}
 \def\bibitem[##1]##2
 {\htmlanchor{##2}{}\RefLabel{##2}[##1]\hangindent=.8cm\hangafter=1}
 }
\def\endreferences{\bigskip\bigskip\endgroup}
\ifundefined{bibstylemodification}\relax\else\bibstylemodification\fi
\startreferences

\bibitem[Billingsley (1995)]{MR1324786}
\bibauthor{Billingsley, P.} (1995).
\newblock {\em Probability and Measure}.
\newblock Wiley Series in Probability and Mathematical Statistics. John Wiley,
  New York, third edition.

\bibitem[Cram\'er and Wold (1936)]{CramerWold}
\bibauthor{Cram\'er, H. \and{} Wold, H.} (1936).
\newblock Some theorems on distribution functions.
\newblock {\em Journal of the London Mathematical Society}, {\bf s1--11}(4),
  290--294.

\bibitem[Crofton (1868)]{Crofton}
\bibauthor{Crofton, M.W.} (1868).
\newblock On the theory of local probability, applied to straight lines drawn
  at random in a plane; the methods used being also extended to the proof of
  certain new theorems in the integral calculus.
\newblock {\em Philos. Trans. Royal Soc. London}, {\bf 158}, 181--199.

\bibitem[John (1955)]{MR0075429}
\bibauthor{John, F.} (1955).
\newblock {\em Plane Waves and Spherical Means Applied to Partial Differential
  Equations}.
\newblock Interscience Publishers, New York-London.

\bibitem[Lyons (2013)]{Lyons:dcov}
\bibauthor{Lyons, R.} (2013).
\newblock Distance covariance in metric spaces.
\newblock {\em Ann. Probab.}, {\bf 41}(5), 3284--3305.

\bibitem[Pollard (2002)]{MR1873379}
\bibauthor{Pollard, D.} (2002).
\newblock {\em A User's Guide to Measure Theoretic Probability}, volume 8 of
  {\em Cambridge Series in Statistical and Probabilistic Mathematics}.
\newblock Cambridge University Press, Cambridge.

\bibitem[Radon (1917)]{Radon}
\bibauthor{Radon, J.} (1917).
\newblock {\"U}ber die {B}estimmung von {F}unktionen durch ihre {I}ntegralwerte
  l\"angs gewisser {M}annigfaltigkeiten.
\newblock {\em Ber. Verh. S\"achs. Akad. Wiss. Leipzig, Math. Nat. kl.}, {\bf
  69}, 262--277.

\bibitem[R{\'e}nyi (1952)]{MR0053422}
\bibauthor{R{\'e}nyi, A.} (1952).
\newblock On projections of probability distributions.
\newblock {\em Acta Math. Sci. Hungar.}, {\bf 3}, 131--142.

\bibitem[Schneider and Weil (2008)]{MR2455326}
\bibauthor{Schneider, R. \and{} Weil, W.} (2008).
\newblock {\em Stochastic and Integral Geometry}.
\newblock Probability and its Applications (New York). Springer-Verlag, Berlin.

\bibitem[Walther (1997)]{walther:withAdd}
\bibauthor{Walther, G.} (1997).
\newblock On a conjecture concerning a theorem of {C}ram\'er and {W}old.
\newblock {\em J. Multivariate Anal.}, {\bf 63}(2), 313--319.
\newblock Addendum, {\it J. Multivariate Anal.} {\bf 67}(2) (1998), 431.

\endreferences